\documentclass[12pt]{amsart}
\usepackage{epsfig,color}
\usepackage{amssymb}

\makeatother

\theoremstyle{definition}

\def\fnum{equation}

\numberwithin{equation}{section}

\newcommand{\Vol}{{\text{Vol}}}

\newcommand{\nn}{{\bf{n}}}

\newcommand{\Hess}{{\text {Hess}}}

\def\RR{{\bold R}}

\newcommand{\cS}{{\mathcal{S}}}

\title[Level set method]{Level set method\\
  \footnotesize\mdseries
  \tiny for motion by mean curvature
 }

\author[]{Tobias Holck Colding}%
\author[]{William P. Minicozzi II}%
 
\thanks{Tobias Holck Colding is Cecil and Ida B. Green Distinguished professor of mathematics at MIT. His email address is  colding@math.mit.edu.
William P. Minicozzi II is professor of mathematics at MIT. His email address is minicozz@math.mit.edu.}

\begin{document}

\maketitle

\begin{abstract}
Modeling of a wide class of physical phenomena, such as crystal growth and flame propagation,
 leads to tracking fronts moving with curvature-dependent speed.  When the speed is the curvature this leads to one of the classical degenerate
nonlinear second order differential equations on Euclidean space. One naturally wonders ``what is the regularity of solutions?''  A priori solutions are only
defined in a weak sense, but it turns out that they are always twice differentiable classical
solutions. This result is optimal;   
  their second derivative is continuous only  in very rigid situations that have a simple geometric interpretation.
The proof weaves together analysis and geometry. Without deeply understanding
the underlying geometry, it is impossible to prove fine analytical properties.
\end{abstract}


\vskip4mm
The spread of  a forest fire, the  growth of a crystal, 
 the inflation of an airbag, 
 and a droplet of oil floating in water  can all be modeled by the level set method.  One of the challenges for modeling is  the presence of discontinuities.
For example, two separate fires can expand  over time and eventually merge to one as in Figure \ref{f:figure1} or a droplet of liquid can split     as in Figure \ref{f:figure1a}.  The level set method allows for this.
 
 \begin{figure}
       \centering\includegraphics[totalheight=.35\textheight, width=.9\textwidth]{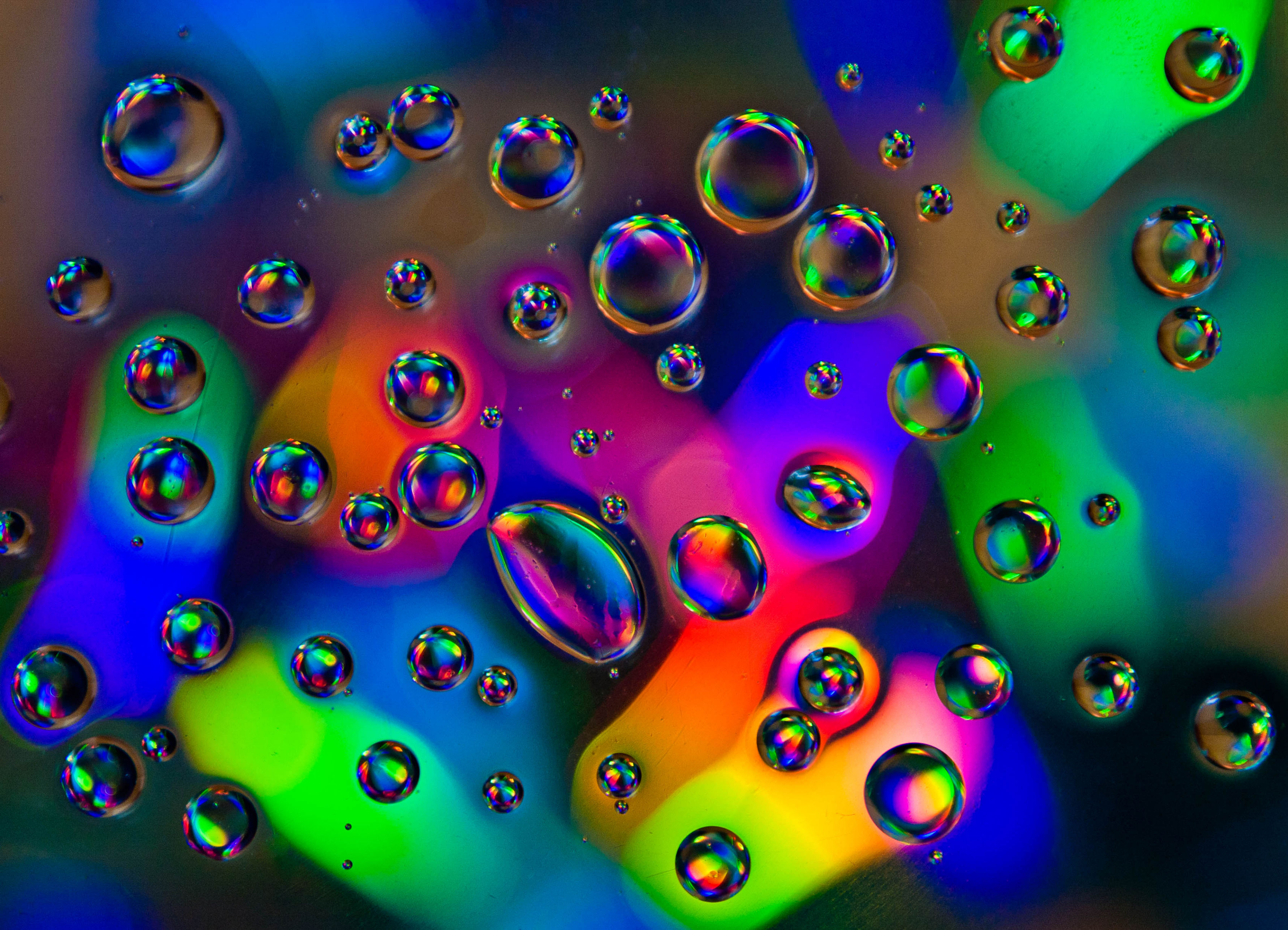}
   \caption{Oil droplets in water can be modeled by the level set method;   https://www.flickr.com/photos/ltn100/7281662212.}
 \label{f:figure0}
\end{figure}

 \begin{figure}[htbp]
 
    \centering\includegraphics[totalheight=.2\textheight, width=.9\textwidth]{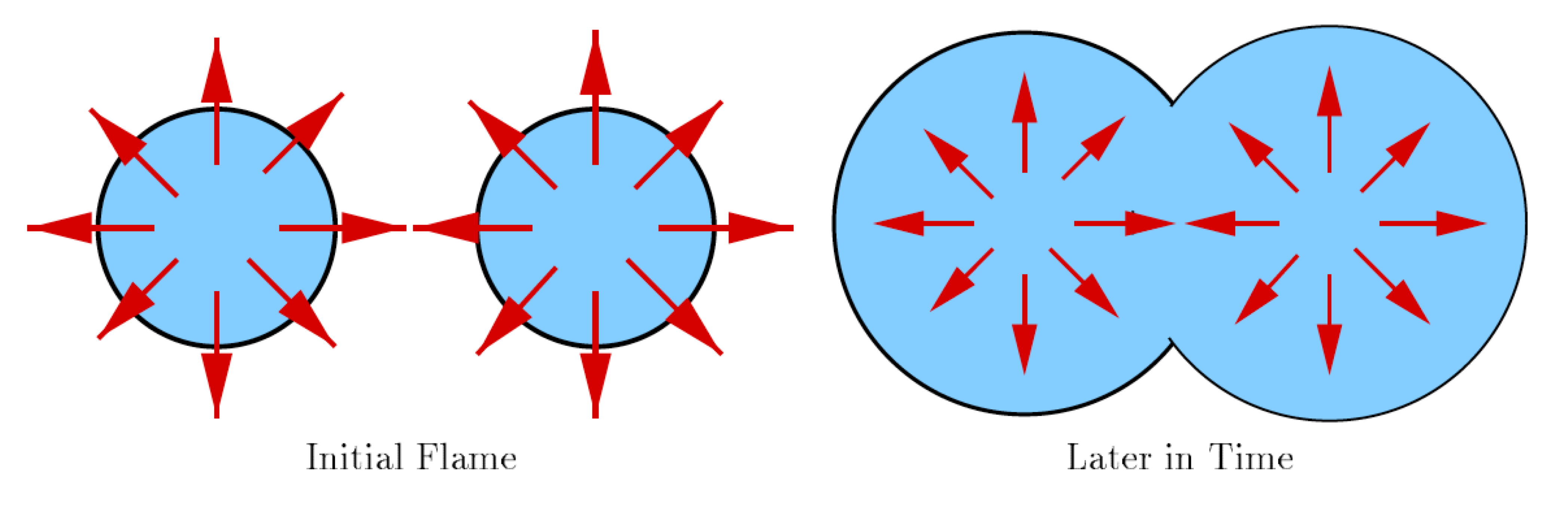}
   \caption{After two fires merge, the evolving front is connected. Figure created by James Sethian and used with permission.}
 \label{f:figure1}
\end{figure}

 \begin{figure}[htbp]
 
    \centering\includegraphics[totalheight=.3\textheight, width=.3\textwidth]{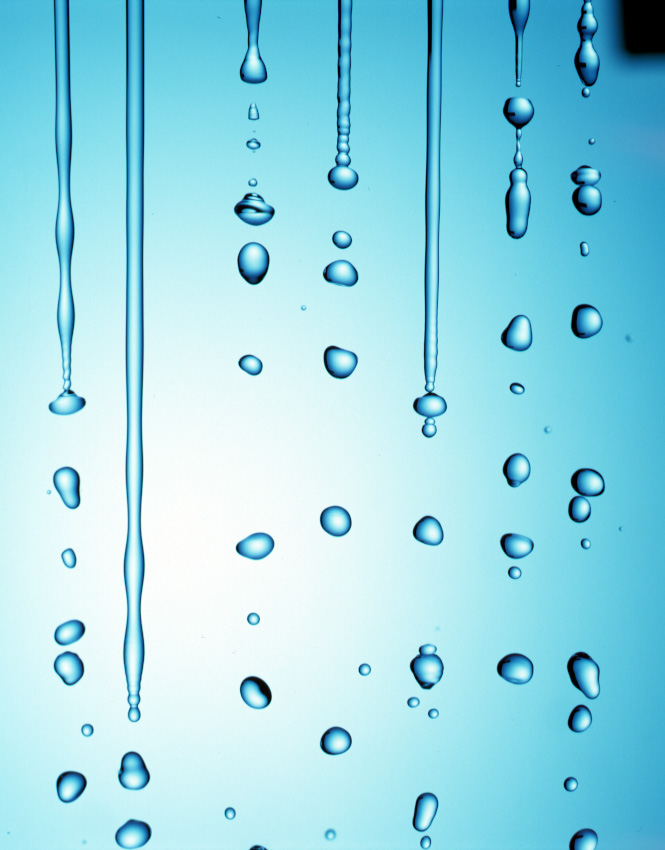}
   \caption{Water drops;  Adam Hart-Davis, http://www.adam-hart-davis.org.}
 \label{f:figure1a}
\end{figure}

\begin{figure}[htbp]
 
   \centering\includegraphics[totalheight=.3\textheight, width=.6\textwidth]{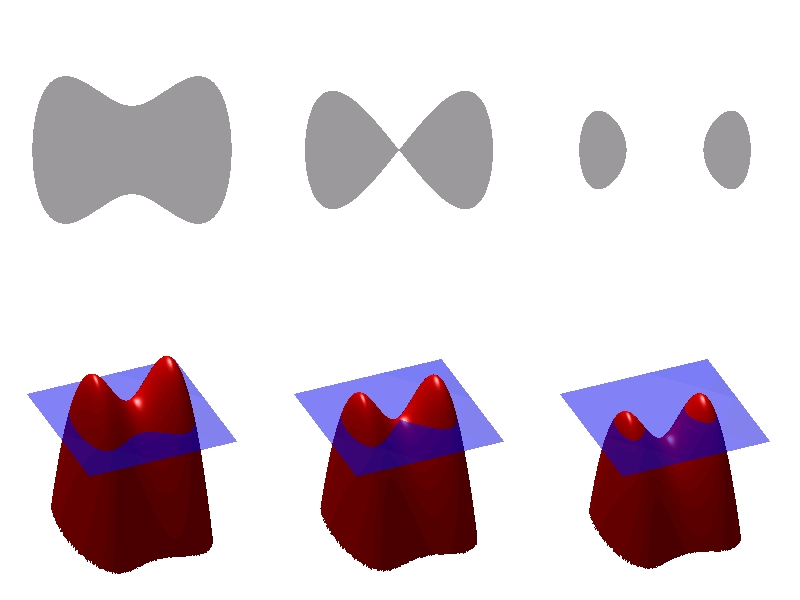}
\caption{
The gray areas represent trees that a forest fire has not yet reached. The edge is the burning fire front that is moving inward  as time propagates left to   right in the upper part of the figure.  The propagation of the fire front is given as a level set of an evolving function in the second line.  Figure from
https://commons.wikimedia.org/w/index.php?curid=2188899.}
\label{f:fig4}
\end{figure}

 The level set method has been used with great success over the last thirty years in both pure and applied mathematics.
  Given  an initial interface or front $M_0$ bounding a region   in $\RR^{n+1}$, the level set method is used to analyze  its subsequent motion under a velocity field.     The idea is to represent the evolving front as a level set of a function $v(x,t)$, where $x$ is in  $\RR^{n+1}$ and $t$ is time.  
  The initial front $M_0$ is given by $$M_0 = \{ x \, | \, v(x,0) = 0\}\, ,$$  
and the evolving front is described for all later time $t$ as the set where $v(x,t)$ vanishes as in Figure \ref{f:fig4}.
There are many functions that have $M_0$ as a level set, but the evolution of the level set does not depend on the choice of the function $v(x,0)$.
 
 

\begin{figure}[htbp]
 
 \centering\includegraphics[totalheight=.1\textheight, width=.6\textwidth]{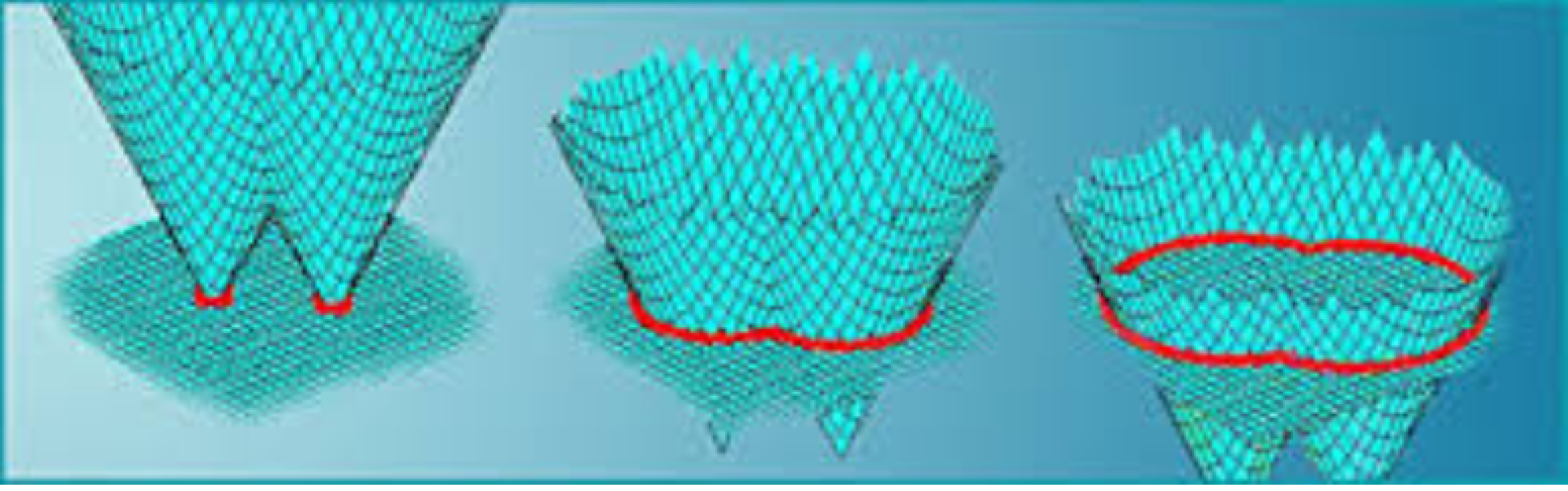}
\caption{The red curves represent  the evolving front at three different times.  Each is the zero level set of the evolving function.
The front here is moving outward as in Figure \ref{f:figure1}.  In Figures \ref{f:figure1a} and \ref{f:fig4}, the front is moving inward. Figure created by James Sethian and used with permission.}
\end{figure}

In mean curvature flow, the velocity vector field is the mean curvature vector and the evolving  front  is  the level set of a function that satisfies a nonlinear degenerate parabolic equation. Solutions are defined in a weak, so-called ``viscosity'' sense; in general, they may not even be differentiable (let alone twice differentiable).
However, it turns out that for a monotonically advancing front viscosity solutions are in fact twice differentiable and satisfy the equation in the classical sense. Moreover, the situation becomes very rigid when the second derivative is continuous.

Suppose $\Sigma\subset \RR^{n+1}$ is an embedded hypersurface and $\nn$ is the unit normal of $\Sigma$.  The {\bf{mean curvature}} is given by $H=\text{div}_{\Sigma}(\nn)$. 
Here $$\text{div}_{\Sigma}(\nn)=\sum_{i=1}^n\langle \nabla_{e_i}\nn,e_i\rangle\ ;$$ where $e_i$ is an orthonormal basis for the tangent space of $\Sigma$.   
  For example, at a point where $\nn$ points in the $x_{n+1}$ direction and the principal directions are in the other axis directions,
\begin{align}
	\text{div}_{\Sigma}(\nn) = \Sigma_{i=1}^n \frac{ \partial \nn_i}{\partial{x_i}} \notag
\end{align}
is the sum ($n$ times the mean) of the principal curvatures.  
If $\Sigma=u^{-1}(s)$ is the level set of a function $u$ on $\RR^{n+1}$ and $s$ is a regular value, then $\nn=\frac{\nabla u}{|\nabla u|}$ and 
$$H=  \sum_{i=1}^n\langle \nabla_{e_i}\nn,e_i\rangle = 
\text{div}_{\RR^{n+1}}\,\left(\frac{\nabla u}{|\nabla u|}\right) \, . $$
 The last equality used that $\langle \nabla_{\nn} \nn , \nn \rangle$ is automatically $0$ because $\nn$ is a unit vector.

A  one-parameter family of smooth hypersurfaces $M_t\subset \RR^{n+1}$ flows by the  {\bf{mean curvature flow}}  if the speed is  equal to the mean curvature and points inward:
\begin{equation}
x_t= -H\,\nn\, , \notag
\end{equation}
where $H$ and $\nn$ are   the mean curvature and unit normal of 
$M_t$ at the point $x$.   Our flows will always start at a smooth embedded connected hypersurface, even if it becomes disconnected and non-smooth at later times.
 The earliest reference to the mean curvature flow we know of is in the work of George Birkhoff in the 1910s, where he used a discrete version of this, and independently in the material science literature in the 1920s.  

 \vskip2mm
  \noindent
  {\bf{Two key properties}}: 
  \begin{itemize}
  \item $H$ is the gradient of area, so mean curvature flow is the negative gradient flow for volume
  ($ \Vol \,M_t$ decreases most efficiently).
  \item  Avoidance property:  If $M_0$ and $N_0$ are disjoint, then $M_t$ and $N_t$ remain disjoint.
  \end{itemize}
  
\begin{figure}[htbp]
\includegraphics[totalheight=.15\textheight, width=.5\textwidth]{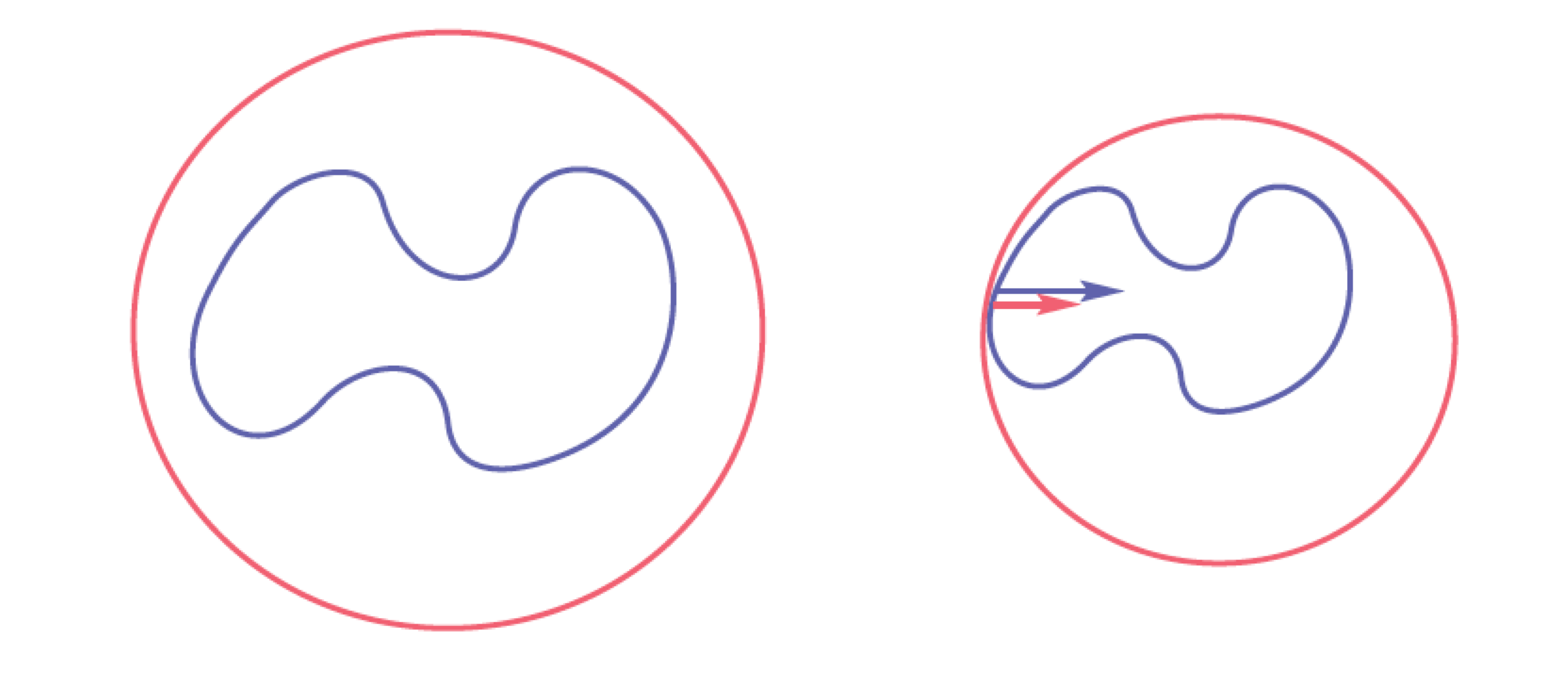}
\caption{Disjoint surfaces avoid each other; contact leads to a contradiction.}
\label{f:avoid}
 \end{figure}

  The avoidance principle is simply a geometric formulation of the maximum principle.  An application of it is illustrated in Figure \ref{f:avoid}, which shows that if one closed hypersurface (the red one) encloses another (the blue one), then the outer one can never catch up with the inner.  The reason for this is that if it did there would be a first point of contact and right before that the inner one would contract faster than the outer, contradicting that the outer was catching up.

 \vskip2mm
 \noindent 
 {\bf{Curve shortening flow}}:  
 When $n=1$ and the hypersurface is a curve, the flow is the curve shortening flow. 
  Under the curve shortening flow, a round circle shrinks through round circles to a point in finite time.
  A remarkable result of  Matthew Grayson from 1987 (using earlier work of Richard Hamilton and Michael Gage) shows that any simple closed curve in the plane remains smooth under the flow until it disappears in finite time in a point.  Right before it disappears,
   the curve will be an almost round circle.    
 
 The evolution of the snake-like simple closed curve  in  Figure \ref{f:spiral}  illustrates this remarkable fact.   (The eight figures are time shots of the evolution.)


\begin{figure}[htbp]
    \begin{minipage}[t]{0.22\textwidth}
    \centering\includegraphics[totalheight=.12\textheight, width=.78\textwidth]{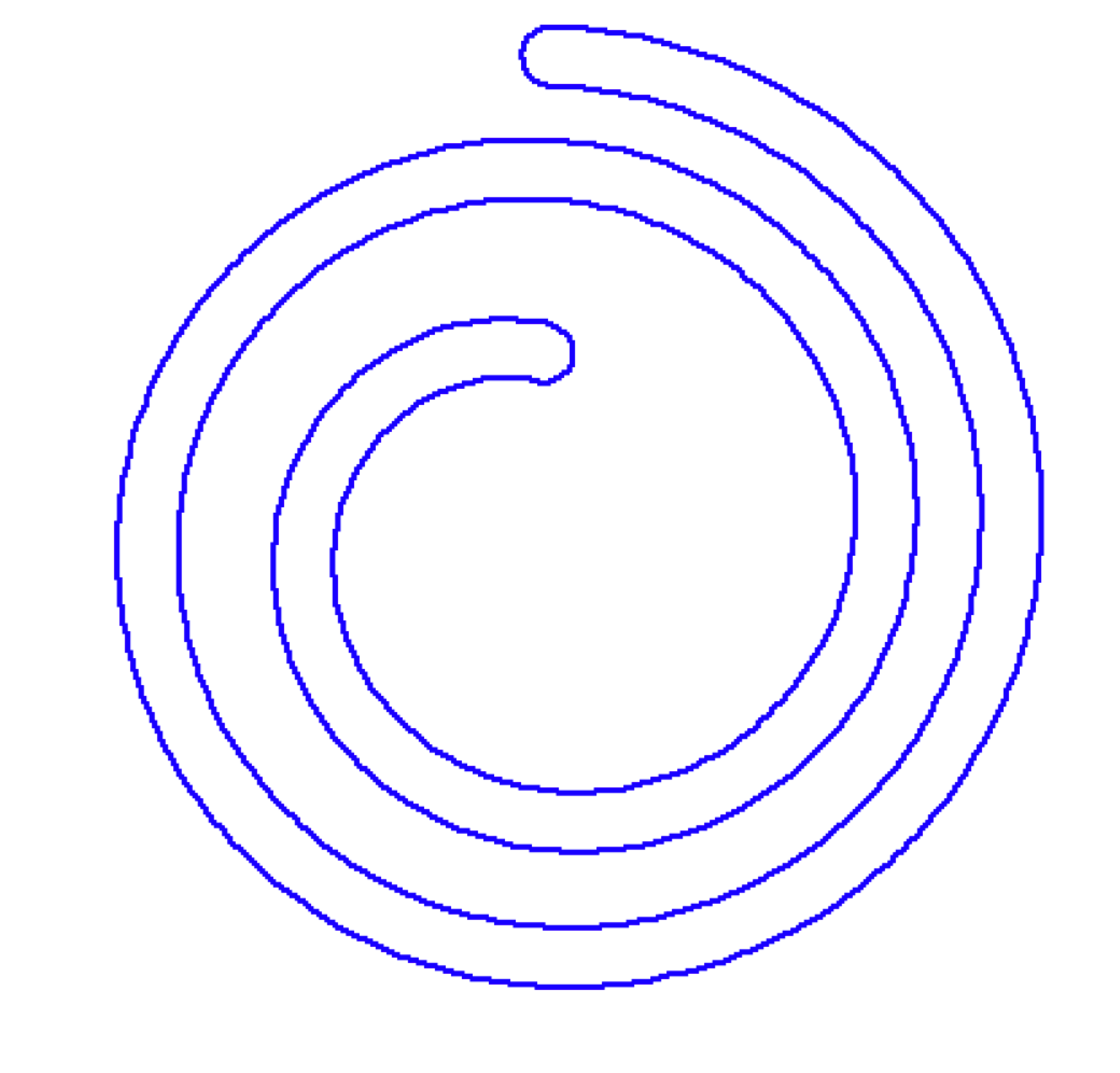} 
    \end{minipage}
    \begin{minipage}[t]{0.22\textwidth}
    \centering\includegraphics[totalheight=.12\textheight, width=.78\textwidth]{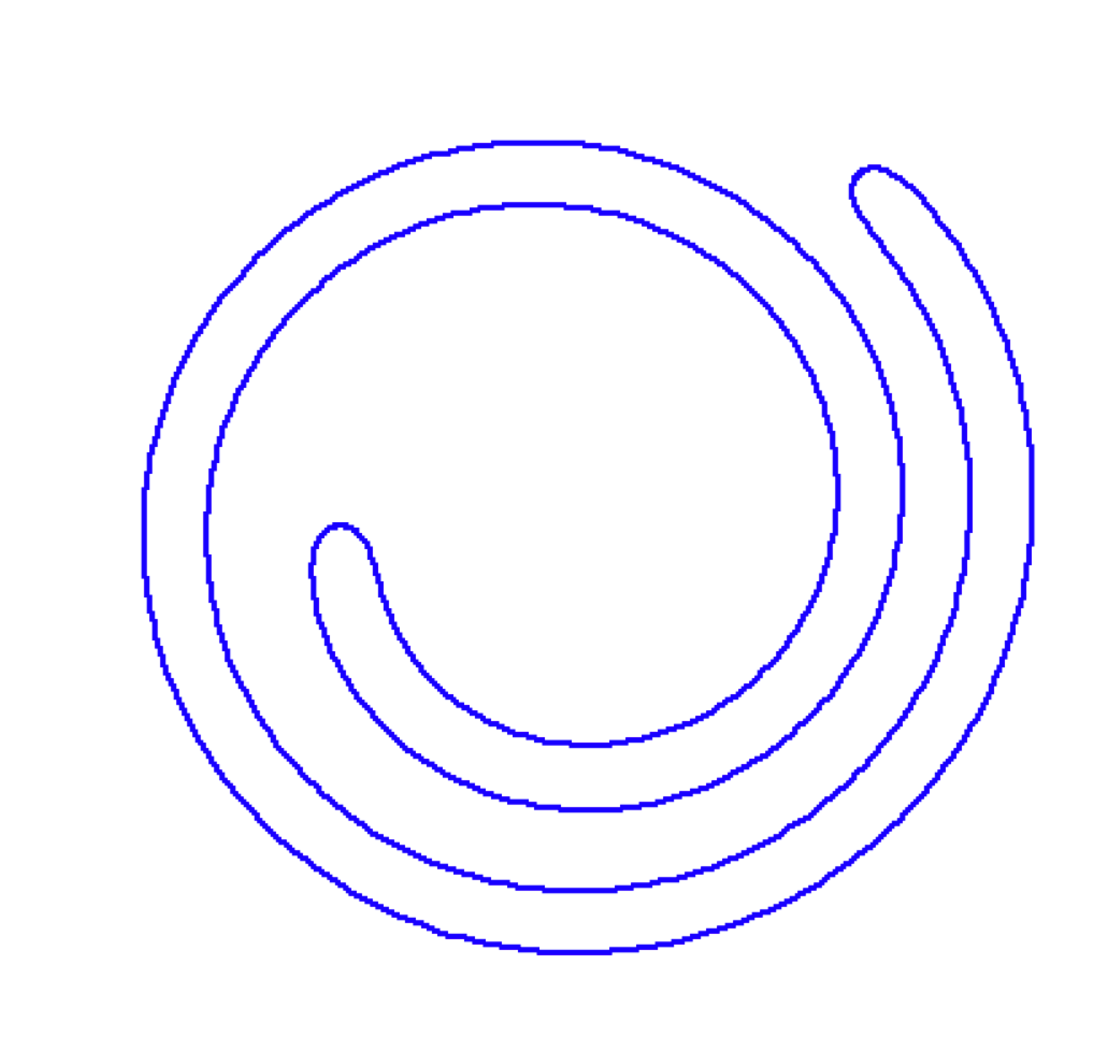}
    \end{minipage}
    \begin{minipage}[t]{0.22\textwidth}
    \centering\includegraphics[totalheight=.12\textheight, width=.78\textwidth]{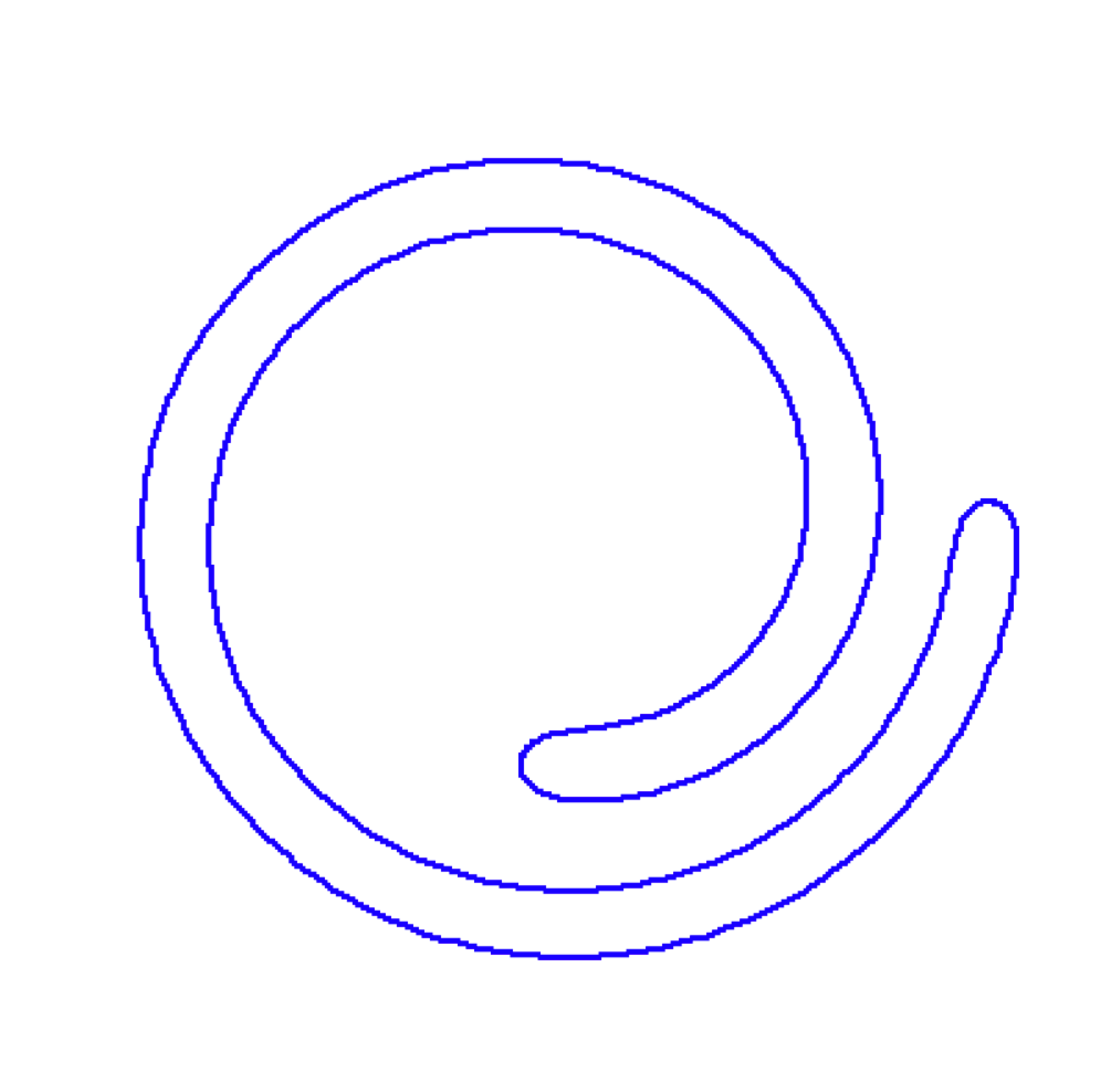}
    \end{minipage}
    \begin{minipage}[t]{0.22\textwidth}
    \centering\includegraphics[totalheight=.12\textheight, width=.78\textwidth]{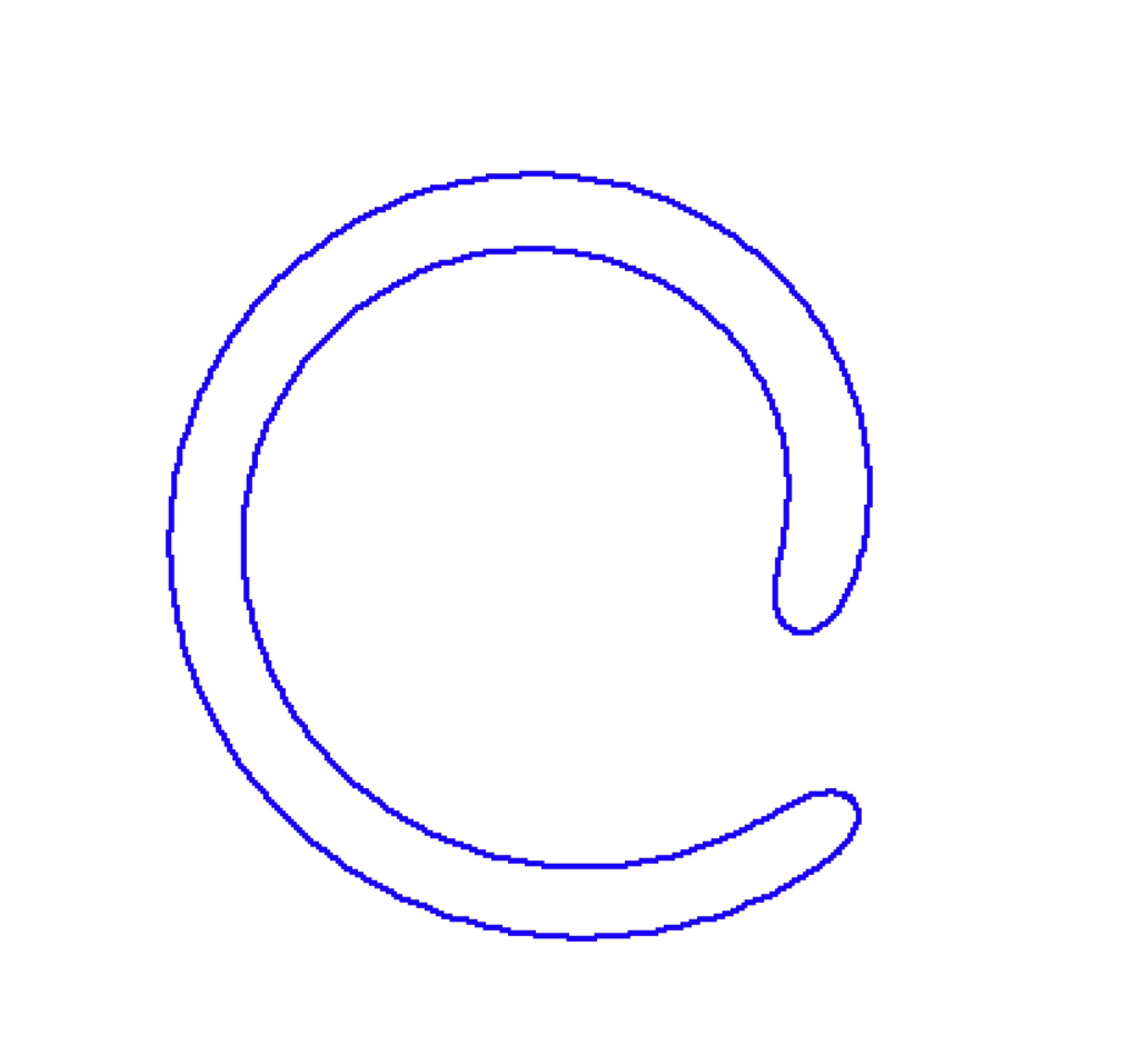} 
    \end{minipage} 
\end{figure}

\begin{figure}[htbp]
    \begin{minipage}[t]{0.22\textwidth}
    \centering\includegraphics[totalheight=.12\textheight, width=.78\textwidth]{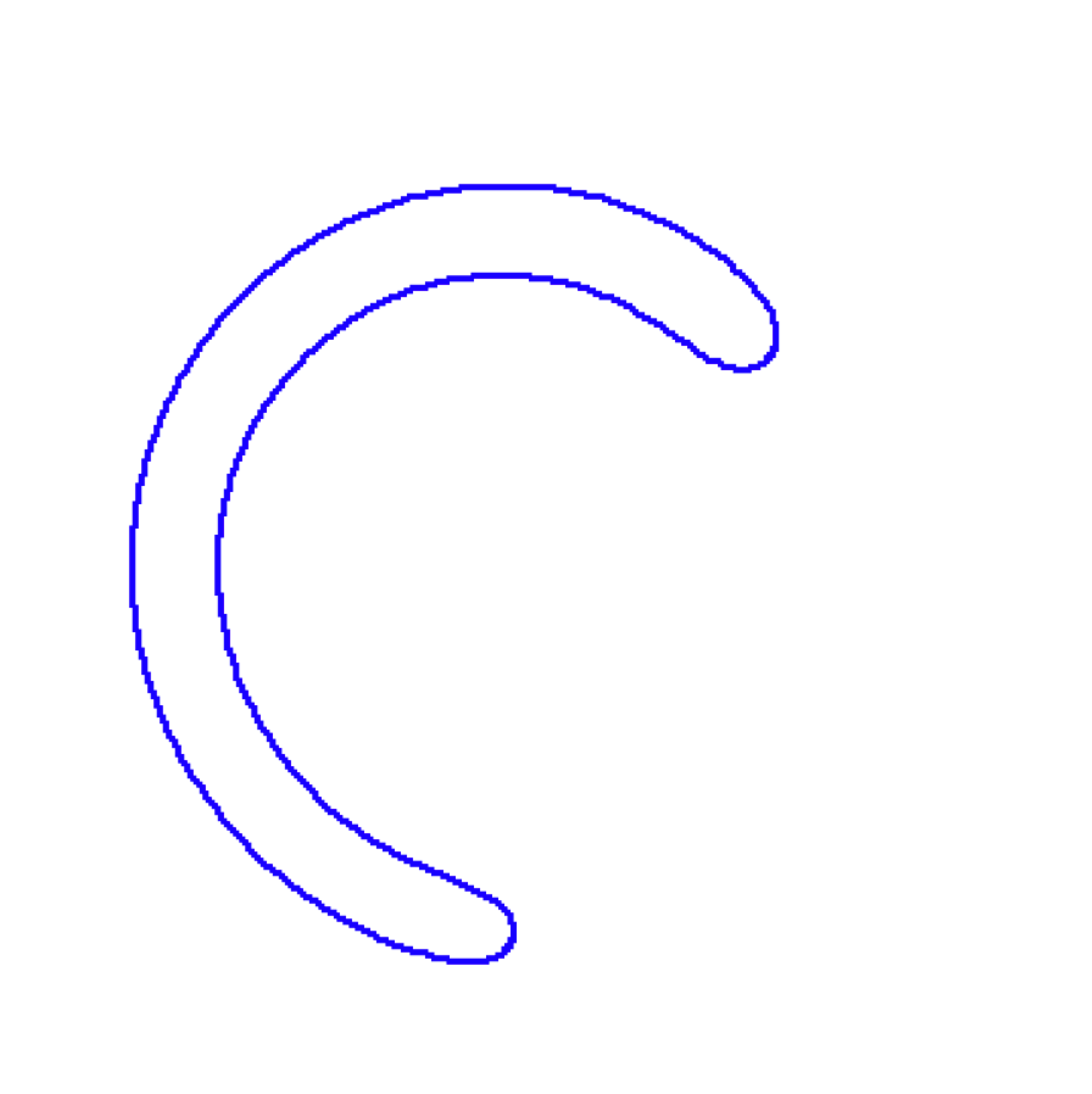} 
    \end{minipage}
    \begin{minipage}[t]{0.22\textwidth}
    \centering\includegraphics[totalheight=.12\textheight, width=.78\textwidth]{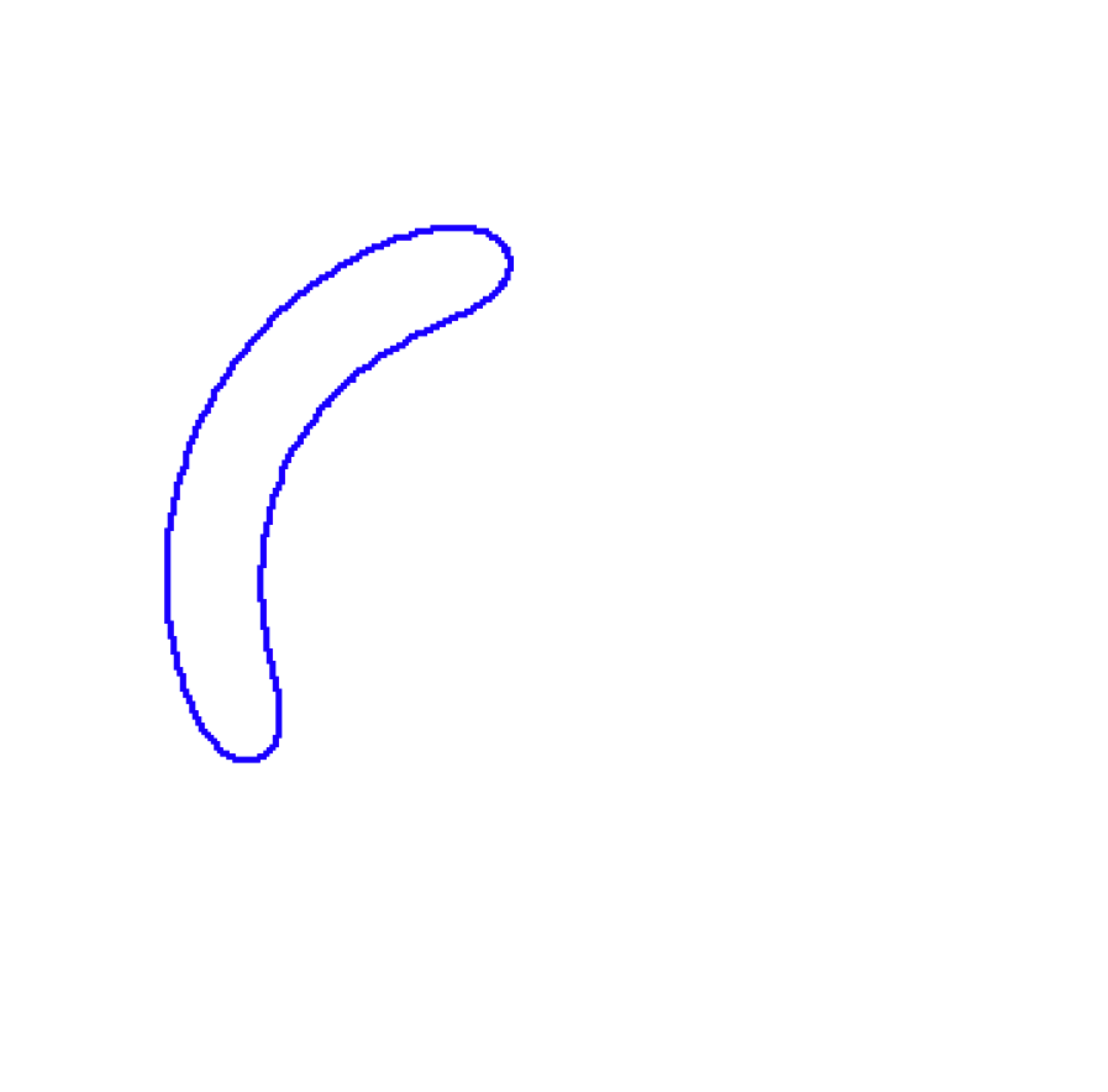}
    \end{minipage}
    \begin{minipage}[t]{0.22\textwidth}
    \centering\includegraphics[totalheight=.12\textheight, width=.78\textwidth]{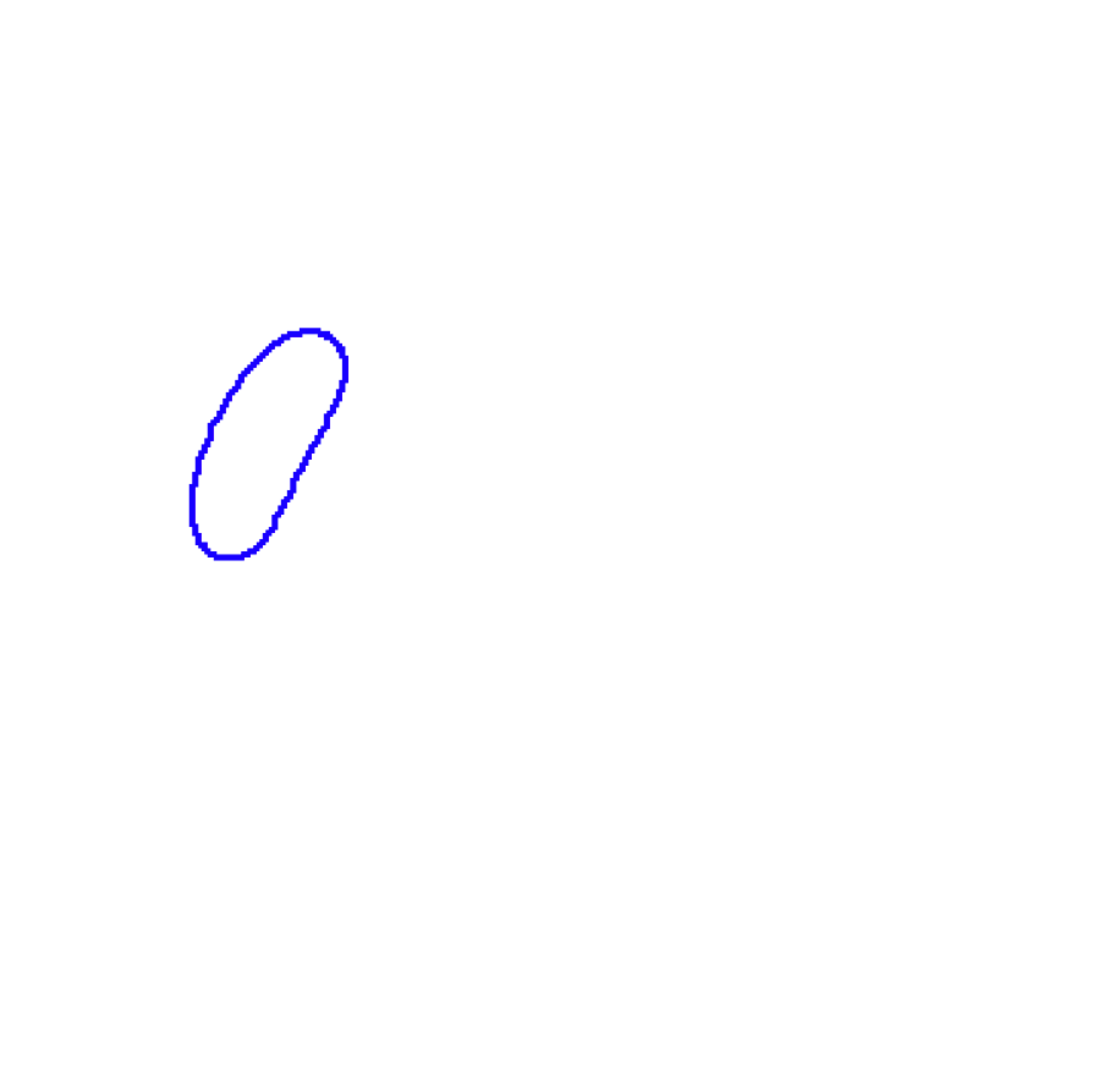}
    \end{minipage}
    \begin{minipage}[t]{0.22\textwidth}
    \centering\includegraphics[totalheight=.12\textheight, width=.78\textwidth]{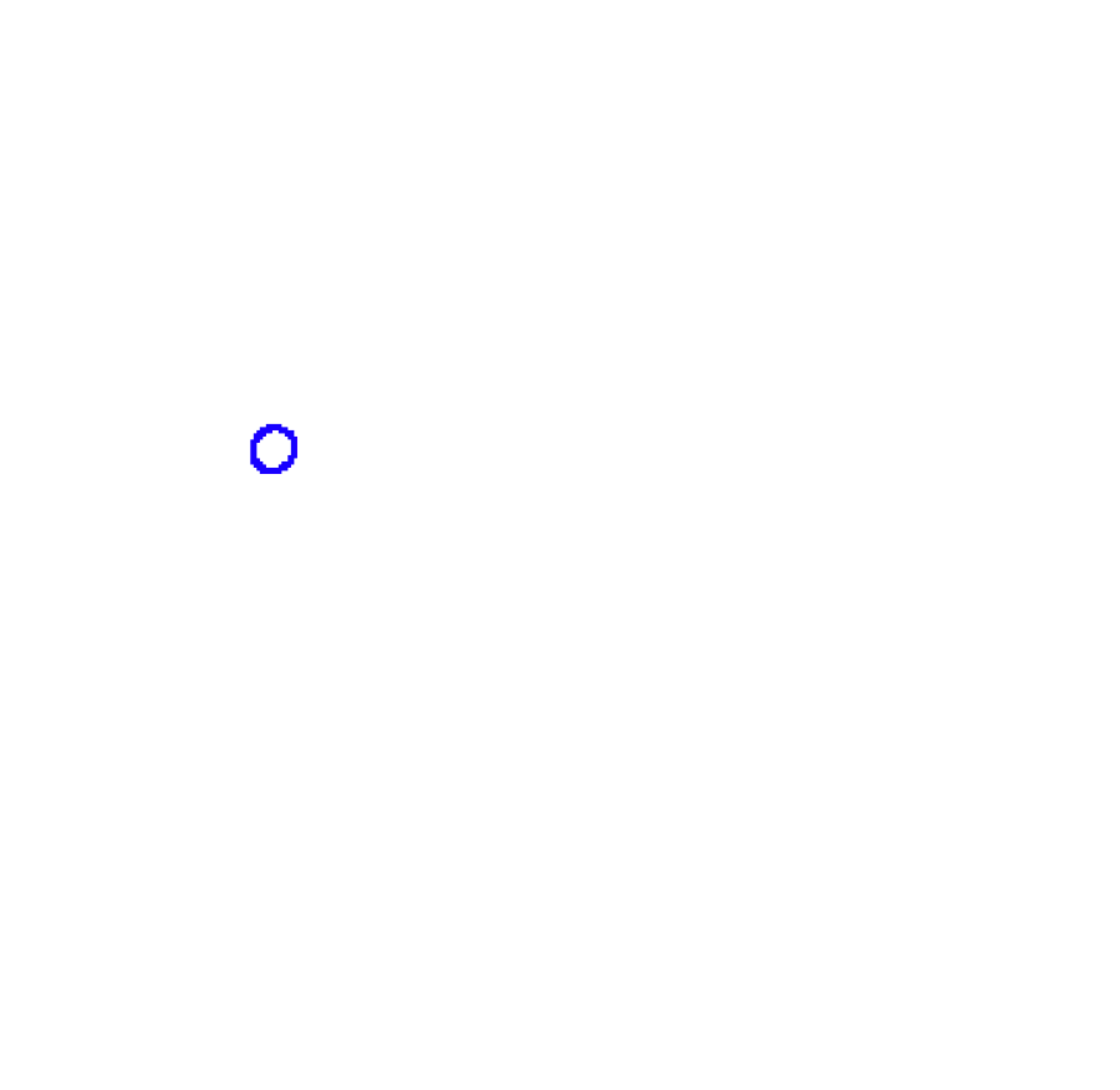} 
    \end{minipage} 
     \caption{Grayson proved that even a tightly wound region becomes round under curve shortening flow. Figures   created by computer simulation by Uwe Mayer  and used with permission.}
     \label{f:spiral}
\end{figure}

\vskip2mm
 \noindent 
 {\bf{Level set flow}}: 
 The analytical formulation of the flow is the level set equation that can be deduced as follows.
  Given a closed embedded hypersurface $\Sigma \subset \RR^{n+1}$, choose a function $v_0:\RR^{n+1}\to \RR$   that is zero on $\Sigma$, positive inside the domain bounded by $\Sigma$, and negative outside.   (Alternatively, choose a function  that is negative inside and positive outside.)
  \begin{itemize}
  \item If we simultaneously flow $\{ v_0 = s_1 \}$ and   $\{ v_0 = s_2 \}$ for $s_1 \ne s_2$, then avoidance implies they  stay disjoint. 
  \item In the level set flow, we look for  $v:\RR^{n+1}\times [0,\infty)\to \RR$ so that each level set $t \to \{ v(\cdot ,t )=s \}$ flows by mean curvature and $v(\cdot,0)=v_0$.
  \item   If   $\nabla v \ne 0$ and the level sets of $v$ flow by mean curvature, then 
  $$v_t=|\nabla v|\,\text{div}\left(\frac{\nabla v}{|\nabla v|}\right) \, .$$
  \end{itemize}
   This is degenerate parabolic and undefined when $\nabla v=0$.
    It may not have classical solutions.
  
   In a   paper - cited more than $12000$ times - from 1988, Stanley Osher and James Sethian studied this equation  numerically.  
   The analytical foundation was provided by Craig Evans and Joel Spruck in a series of four papers in the early 1990s and, independently and at the same time, by Yun Gang Chen, Yoshikazu Giga, and Shunichi Goto.  
    Both of these two groups constructed (continuous) viscosity solutions and showed uniqueness.  The notion of viscosity solutions  had been developed by Pierre-Louis Lions and Michael G. Crandall  in the early 1980s.  The work of these two groups on the level set flow was one of the significant applications of this theory.

  \vskip2mm
  \noindent 
  {\bf{Examples of singularities}}: 
  Under mean curvature flow a round sphere remains round but shrinks and eventually becomes extinct in a point.    A round cylinder remains round and eventually becomes extinct in a line.   The marriage ring is the example of a thin torus of revolution in $\RR^3$.   
 Under the flow the marriage ring shrinks to a circle then disappears. 
 
  \begin{figure}[htbp]
\includegraphics[totalheight=.3\textheight, width=.6\textwidth]{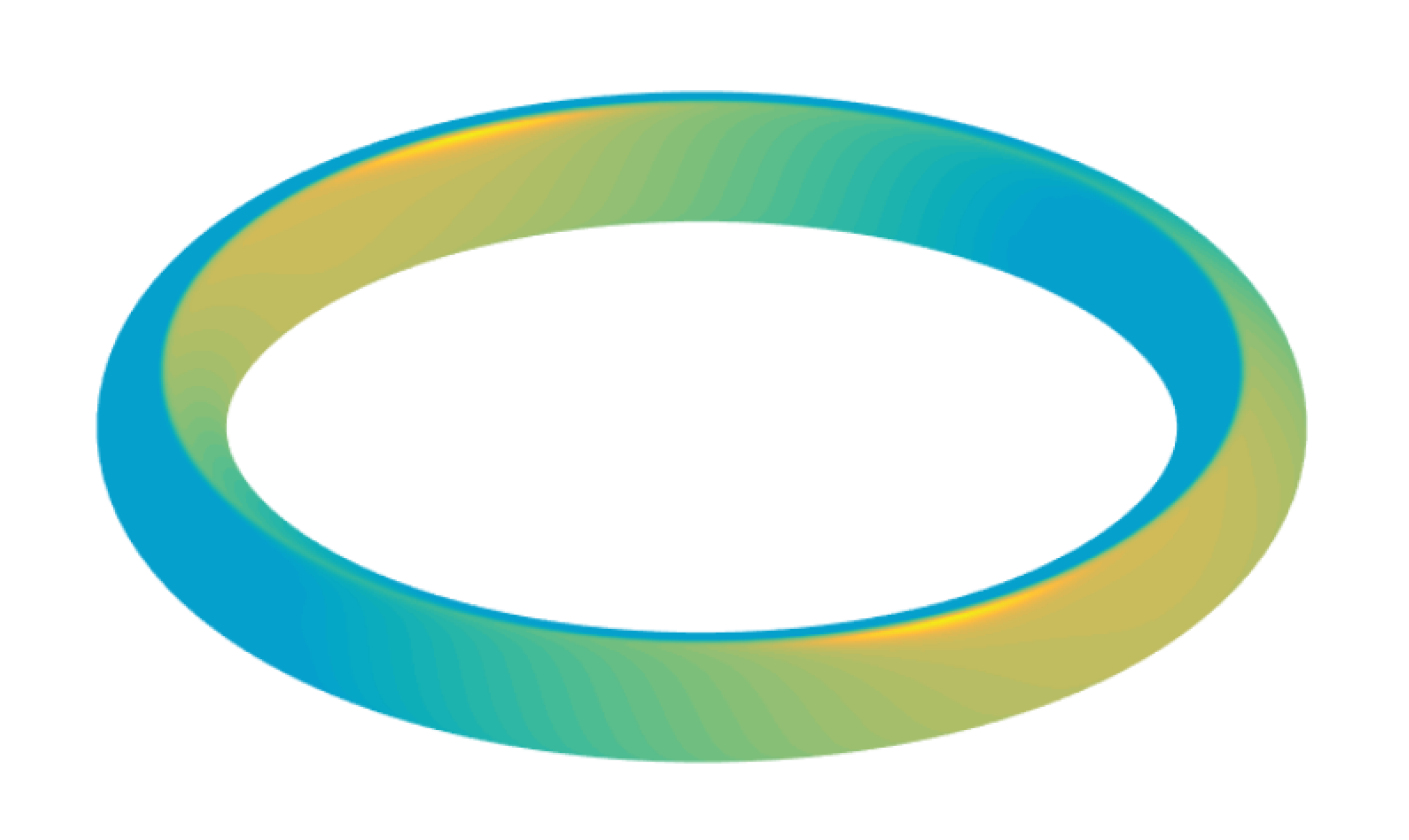}
\caption{The marriage ring shrinks to a circle under mean curvature flow.}
\end{figure}

 \vskip2mm
 \noindent
 {\bf{Dumbbell}}:  
Figure \ref{f:dumb} shows the evolution of a rotationally symmetric mean convex dumbbell in $\RR^3$.  If the neck is sufficiently thin, then the neck pinches off first and the surface disconnects into two components.         Later each component (bell) shrinks to a round point.  This example falls into a larger category of   surfaces that are rotationally symmetric around an axis.  Because of the symmetry, then the solution reduces to a one-dimensional heat equation.   This was analyzed already in the early 1990s by Sigurd Angenent, Steven Altschuler, and Giga; cf. also with work of Halil Mete Soner and Panagiotis Souganidis from around the same time.   A key tool in the arguments of Angenent-Altschuler-Giga was a parabolic Sturm-Liouville theorem of Angenent that holds in one spatial dimension.

\begin{figure}[htbp]
    \begin{minipage}[t]{0.5\textwidth}
    \centering\includegraphics[totalheight=.1\textheight, width=.65\textwidth]{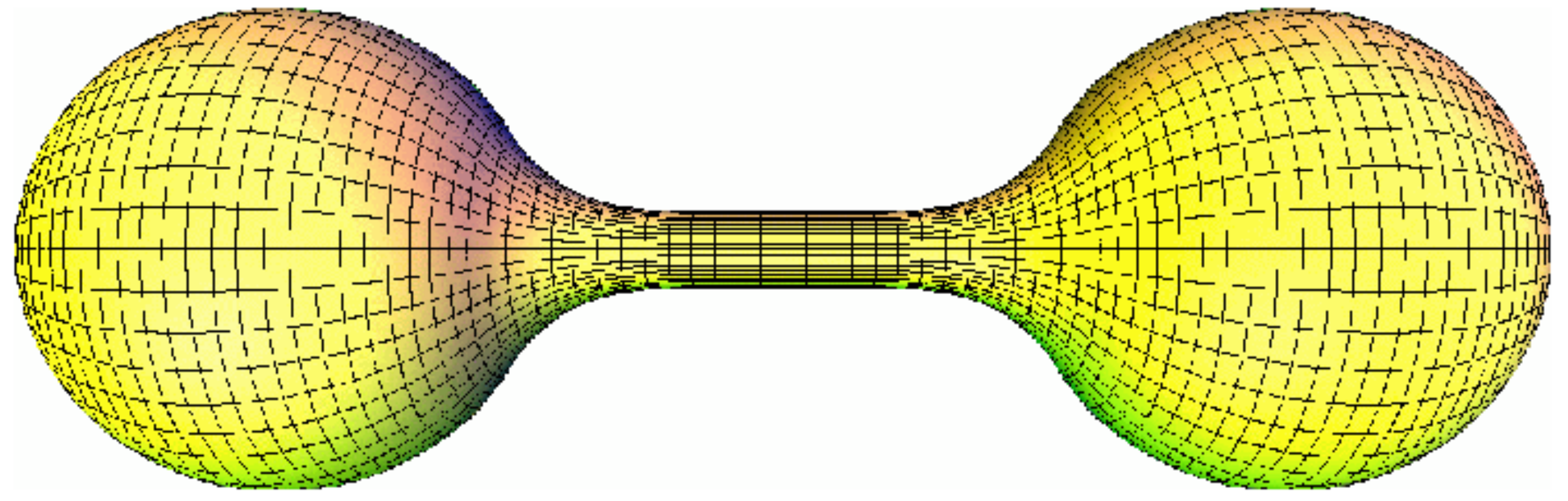}
    \end{minipage}\begin{minipage}[t]{0.5\textwidth}
    \centering\includegraphics[totalheight=.1\textheight, width=.65\textwidth]{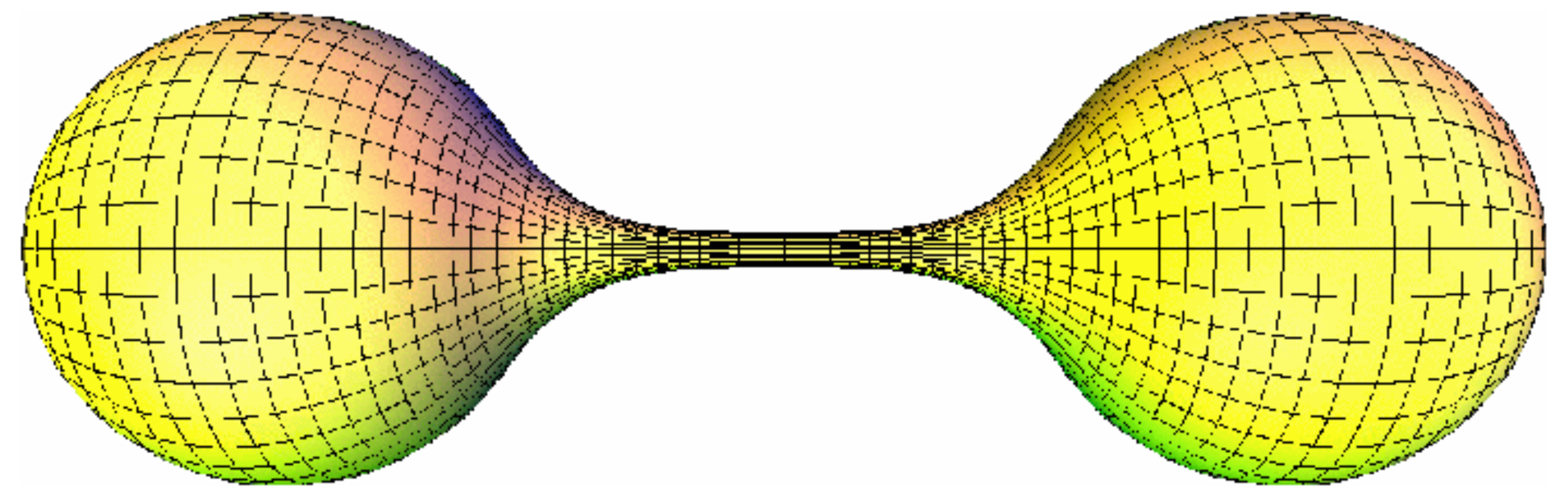}
    \end{minipage}
\end{figure}

\begin{figure}[htbp]
    \begin{minipage}[t]{0.5\textwidth}
    \centering\includegraphics[totalheight=.1\textheight, width=.65\textwidth]{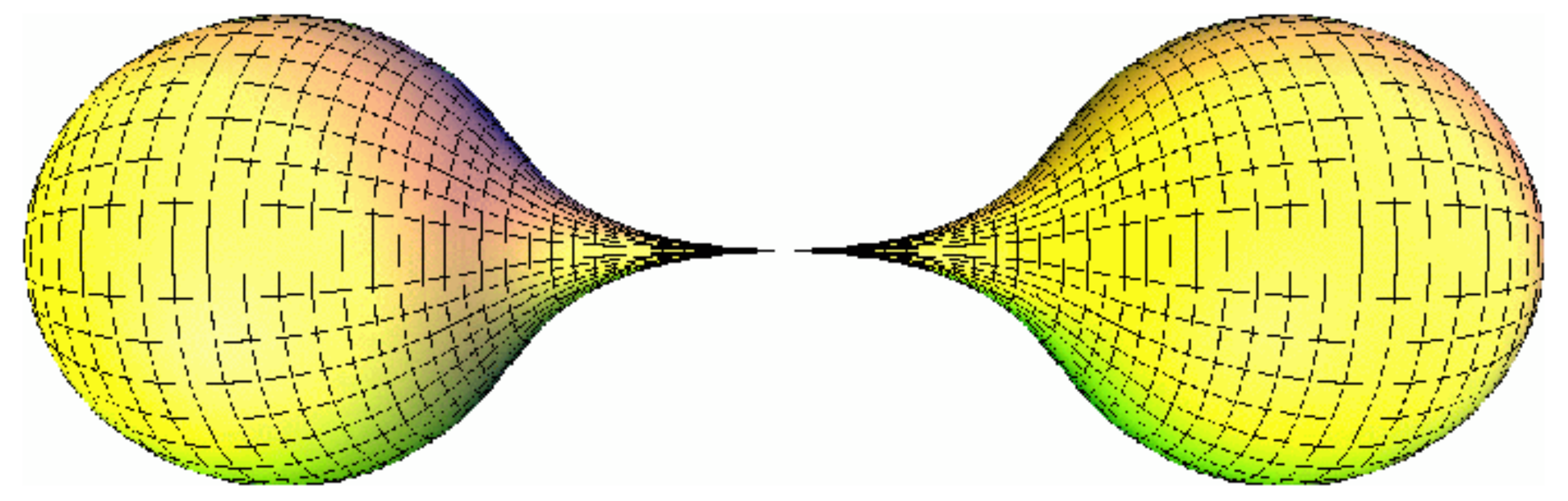}
    \end{minipage}\begin{minipage}[t]{0.5\textwidth}
    \centering\includegraphics[totalheight=.1\textheight, width=.65\textwidth]{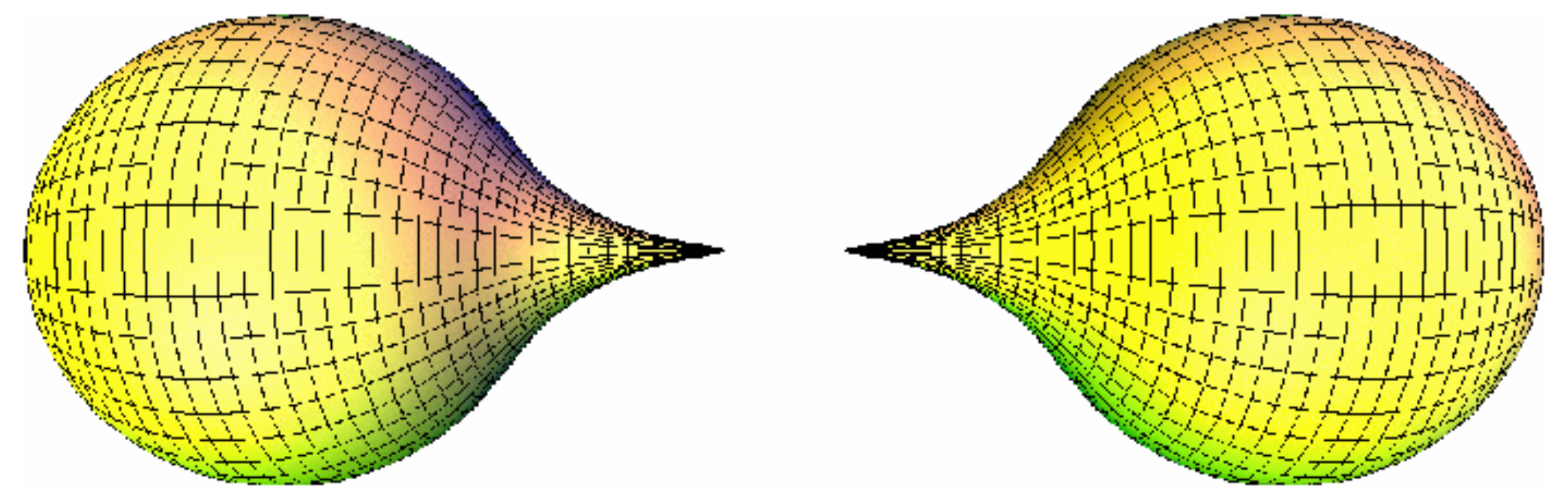}
    \end{minipage}
\end{figure}

\begin{figure}[htbp]
    \begin{minipage}[t]{0.5\textwidth}
    \centering\includegraphics[totalheight=.1\textheight, width=.65\textwidth]{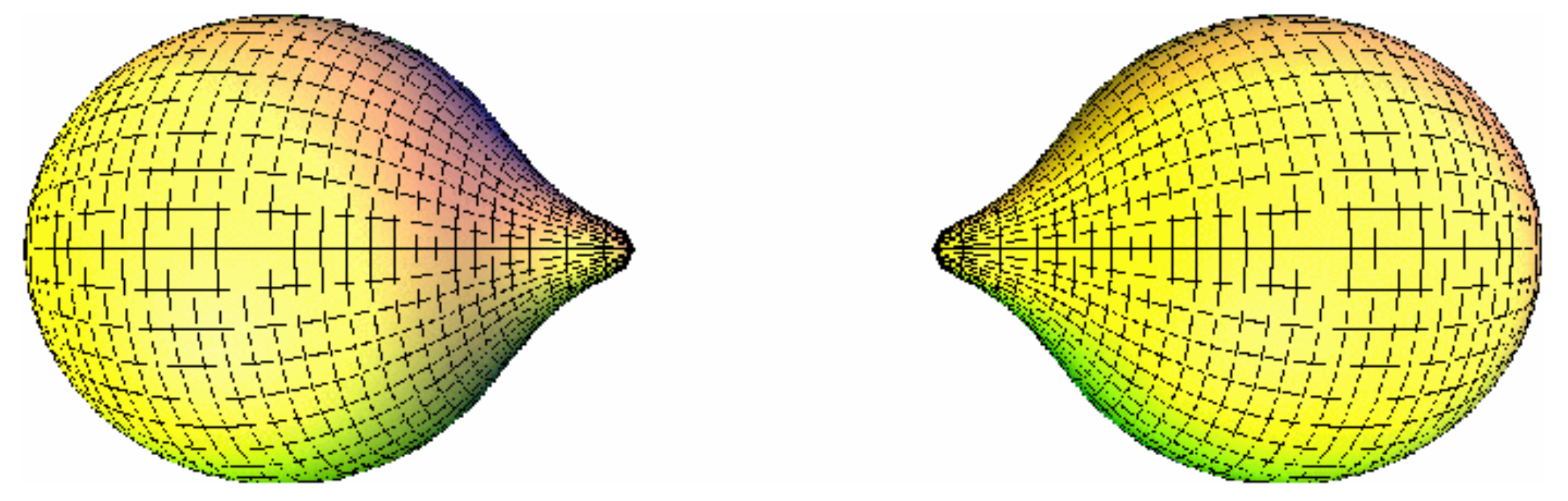}
    \end{minipage}\begin{minipage}[t]{0.5\textwidth}
    \centering\includegraphics[totalheight=.1\textheight, width=.65\textwidth]{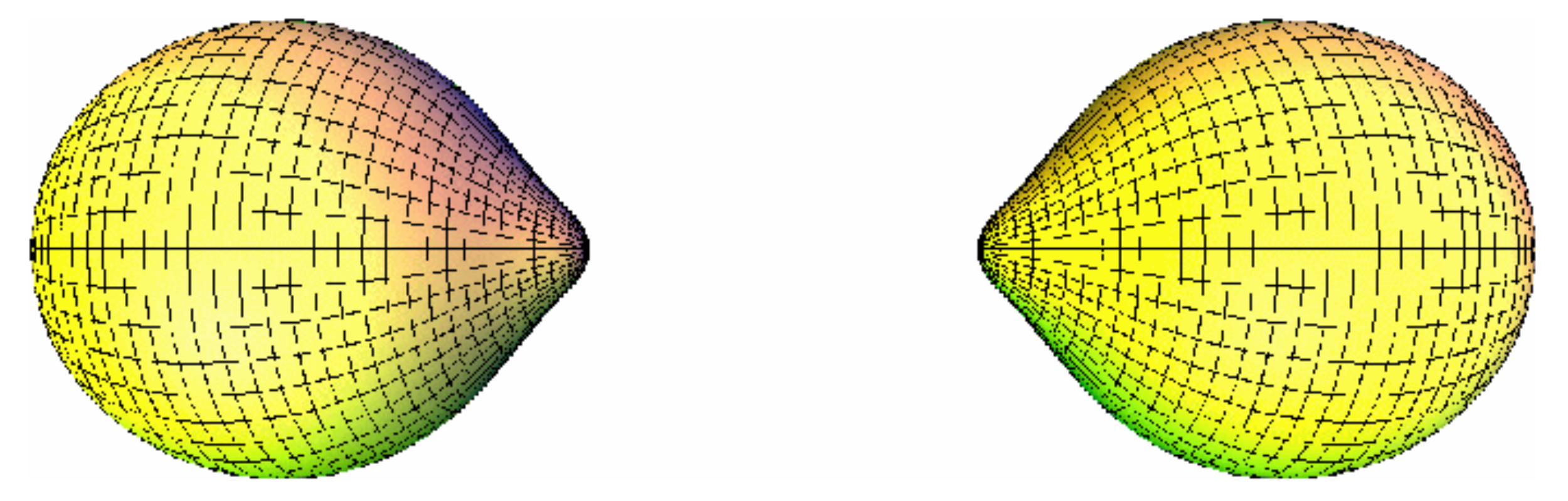}
    \end{minipage}
\end{figure}

\begin{figure}[htbp]
    \begin{minipage}[t]{0.5\textwidth}
    \centering\includegraphics[totalheight=.1\textheight, width=.65\textwidth]{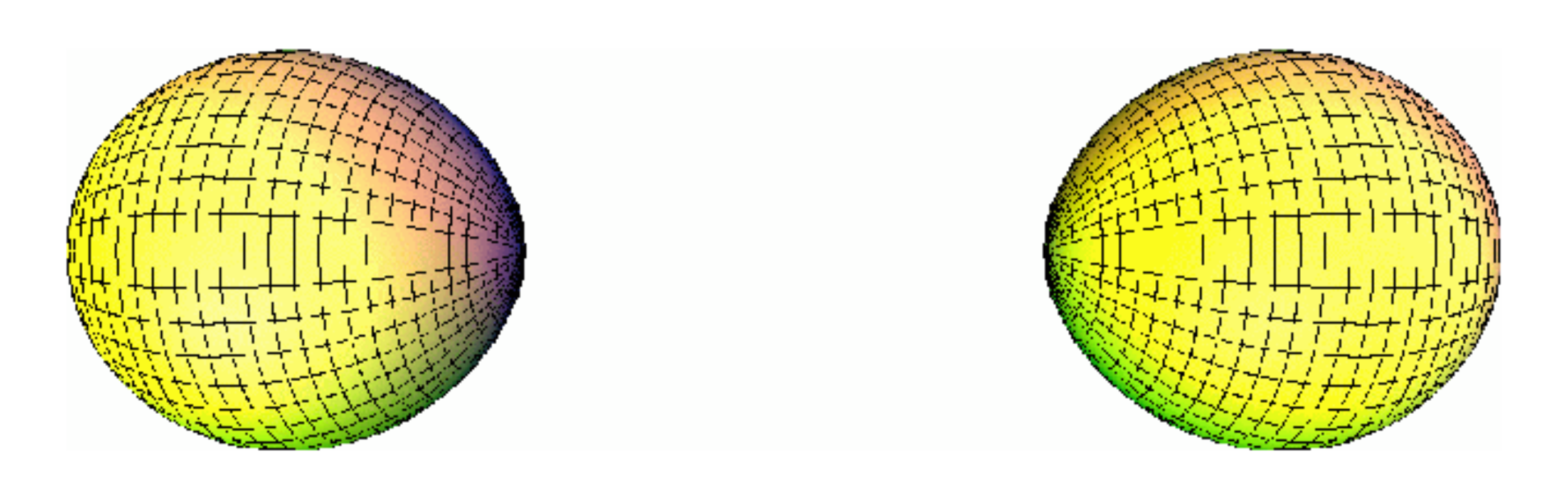}
    \end{minipage}\begin{minipage}[t]{0.5\textwidth}
    \centering\includegraphics[totalheight=.1\textheight, width=.65\textwidth]{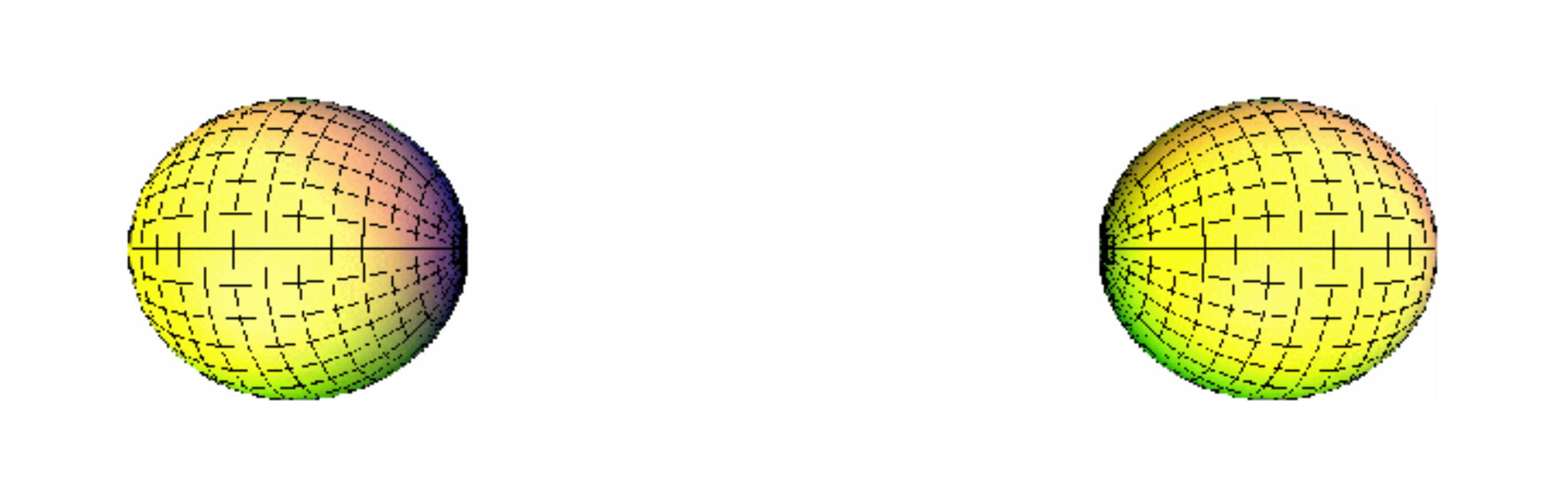}
    \end{minipage}
    \caption{A narrow neck on a dumbbell pinches off, and then each half shrinks to a point. Figures   created by computer simulation by Uwe Mayer  and used with permission.}
     \label{f:dumb}
\end{figure}

  \vskip2mm
  \noindent
  {\bf{Singular set}}:  
    Under mean curvature flow closed hypersurfaces contract, develop singularities and eventually become extinct.   
    The {\bf{singular set}} $\cS$ is the set of points in space and time where the flow is not smooth.
  
  In the first  three examples  - the sphere, the cylinder and the marriage ring -   $\cS$ is  a point, a line, and a closed curve, respectively.  In each case, the singularities occur only at a single time.  
In contrast, the dumbbell has two singular times with one singular point at the first time and two at the second.  
  
 \vskip2mm
\noindent 
 {\bf{Mean convex flows}}: 
 A hypersurface is convex if every principal curvature is positive.  It is mean convex if $H > 0$, i.e., if the sum of the principal curvatures is positive at every point.   
  Under the mean curvature flow,  a mean convex hypersurface moves inward and, since mean convexity is preserved, it will continue to move inward and eventually sweep out the entire compact domain bounded by the initial hypersurface.

  \vskip2mm
  \noindent
  {\bf{Level set flow for mean convex hypersurfaces}}: 
 When the hypersurfaces are mean convex, the equation can be rewritten as a degenerate elliptic equation for a function $u$   defined by
$$u (x) =  \{t\,|\,x\in M_t\} \, .$$
We say that $u$ is the {\bf{arrival time}} since it is  the time the hypersurfaces $M_t$ arrive at $x$ as the front sweeps through the 
compact domain bounded by the initial hypersurface.  The arrival time has a game theoretic interpretation by work of
 Robert Kohn and Sylvia Serfaty.
It follows easily that if we set $v(x,t)=u(x)-t$, then $v$ satisfies the level set flow.  
Now the level set equation $v_t= |\nabla v| \, {\text{div}}(\nabla  v/|\nabla v|)$
 becomes  $$-1=|\nabla u|\,\text{div}\left(\frac{\nabla u}{|\nabla u|}\right)\, .$$
This is a degenerate elliptic equation that is undefined when $\nabla u = 0$.   Note that if $u$ satisfies this equation, then so does $u$ plus a constant.  This just corresponds to shifting the time when the flow arrives by a constant.  
A particular example of a solution to this equation is the function  $u  =-\frac{1}{2}\,(x_1^2+x_{2}^2)$ that is the arrival time for shrinking round cylinders in $\RR^3$.  
In general, Evans-Spruck (cf. Chen-Giga-Goto) constructed Lipschitz solutions to this equation.

  \vskip2mm
  \noindent
  {\bf{Singular set of mean convex level set flow}}:  The singular set of the flow is the critical set of $u$.  Namely, $(x,u(x))$ is singular if and only if $\nabla_xu=0$.
  For instance, in the example of the shrinking round cylinders in $\RR^3$, the arrival time is given by  $u=-\frac{1}{2}\,(x_1^2+x_{2}^2)$  and the flow is singular in the line  $x_1=x_2=0$ that is exactly where $\nabla u=0$.

  \vskip2mm
  We will next see that even though the arrival time was only a solution to the level set equation in a weak sense it turns out to be always twice differentiable classical solution.
  
  \noindent
  {\bf{Differentiability}}: 
\cite{CM2}:
\begin{itemize}
\item $u$ is twice differentiable everywhere, with bounded second derivatives, and smooth away from the critical set.
 \item  $u$ satisfies the equation everywhere in the classical sense.  
 \item  At each critical point the hessian is symmetric and has only two eigenvalues $0$ and $-\frac{1}{k}$;   $-\frac{1}{k}$ has multiplicity $k+1$.    
 \end{itemize}
 
This result is equivalent to saying that at a critical point, say $x=0$ and $u(x) =0$, the function $u$ is (after possibly a rotation of $\RR^{n+1}$) up to higher order terms equal to the quadratic polynomial 
 $$-\frac{1}{k}\,\left(x_1^2+\cdots +x_{k+1}^2\right)\, .$$  This second order approximation is simply the arrival time of the shrinking round cylinders.  It suggests that the level sets of $u$ right before the critical value and near the origin should be approximately cylinders (with an $n-k$ dimensional axis).  This has indeed been known for a long time  and is due to Huisken,  White,   Sinestrari,    Andrews,  Haslhofer-Kleiner.
  It also suggests that those cylinders should be nearly the same (after rescaling to unit size).  That is, the axis of the cylinders should not depend on the value of the level set.  This last property however was only very recently established in \cite{CM1}  and is the key to proving that the function is twice differentiable\footnote{Uniqueness of the axis is parallel to the fact that a function is differentiable at a point precisely if on all sufficiently small scales at that point it looks like the {\it{same}} linear function.}.     The proof  that the axis is unique, independent on the level set, relies on a key new inequality that draws its inspiration from real algebraic geometry although the proof is entirely new.    This kind of uniqueness is a famously difficult problem in geometric analysis and no general case had previously been known.

   \begin{figure}[htbp]		
\centering\includegraphics[totalheight=.3\textheight, width=.3\textwidth]{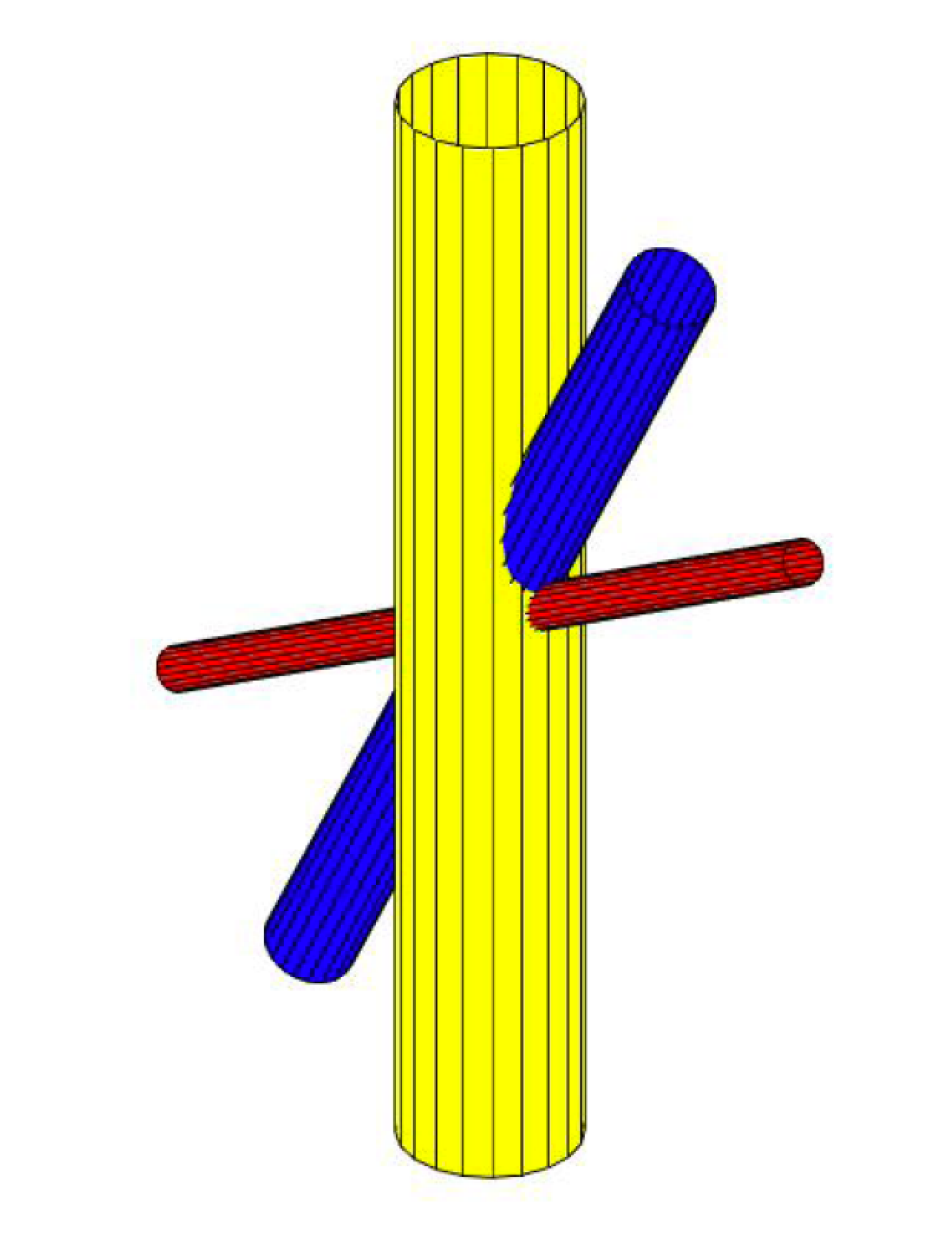}
\caption{The figure illustrates a situation, that turns out to be impossible, where the axes of the cylinders for three different level sets near the critical point are radically different.  
Instead, as the level sets approach the critical point, the axes go to a unique limit.}   \label{f:unique}
  \end{figure}

\vskip2mm
\noindent
{\bf{Regularity of solutions}}:  We have seen that the arrival time is always twice differentiable and one may wonder whether there is even more regularity.  Gerhard Huisken showed already in 1990 
that the arrival time is $C^2$ for {\it{convex}} $M_0$.  However, in 1992 Tom Ilmanen gave an example of a rotationally symmetric {\it{mean convex}} $M_0$ in $\RR^3$ where $u$ is not $C^2$.   This result of Ilmanen shows that the above theorem about differentiability can not be improved to $C^2$.  We will see later that in fact one can entirely characterize when the arrival time is $C^2$.   In the plane, Kohn and Serfaty (2006) showed that $u$ is $C^3$ and for $n>1$ Natasa Sesum   (2008) gave an example of a {\it{convex}} $M_0$ where $u$ is not $C^3$.  Thus Huisken's result is optimal for $n>1$.

\vskip2mm
The next result shows  that one can entirely characterize when the arrival time is $C^2$.

\vskip2mm
\noindent
{\bf{Continuous differentiability}}: \cite{CM3}:  $u$ is $C^2$ if and only if:
\vskip2mm
\begin{itemize}
\item There is exactly one singular time (where the flow becomes extinct).  
\item The singular set $\cS$ is a $k$-dimensional closed connected embedded $C^1$ submanifold of cylindrical singularities.
\end{itemize}
 Moreover, the axis of each cylinder is the tangent plane to $\cS$.

\vskip2mm
When  $u$ is $C^2$ in $\RR^3$,  the singular set $\cS$ is either:
\begin{enumerate}
\item A single point with a spherical singularity.
\item A simple closed $C^1$ curve of  cylindrical singularities.
\end{enumerate}

The examples of the sphere and the marriage ring show that each of these phenomena can happen, whereas the example of the dumbbell does not fall into either, showing that in that case the arrival time is not $C^2$.

\vskip2mm
We can restate this result for $\RR^3$ in terms of the structure of the critical set and Hessian:  $u$ is $C^2$ if and only if $u$ has exactly one critical value  and  the critical set is either:
\begin{enumerate}
\item A single point where $\Hess_u$ is $- \frac{1}{2}$ times the identity.
\item A simple closed $C^1$ curve where $\Hess_u$ has eigenvalues $0$ and $-1$ with multiplicities $1$ and $2$, respectively.
\end{enumerate}
In case (2), the kernel of $\Hess_u$ is tangent to the curve.
 
\vskip2mm
\noindent
{\bf{Concluding remarks}}:  We have seen that for one of the classical differential equations in order to understand the analysis it is necessary to understand the underlying geometry.  There are many tantalizing parallels to other differential equations both elliptic and parabolic. 

\vskip4mm
For references about mean curvature flow, see the survey \cite{CMP} and  \cite{CM1}-\cite{CM3}.

\end{document}